\documentclass[a4paper,12pt]{scrartcl}
\usepackage{mathtext}                    
\usepackage{cmap}                        
\usepackage[cp1251]{inputenc}            
\usepackage[english]{babel}     
\usepackage[left=3cm,right=2cm,top=2cm,bottom=3cm]{geometry} 
\usepackage{amssymb}
\usepackage{amsmath}

\begin{document}

\
\centerline{ON SUMS AND PRODUCTS OF PERIODIC FUNCTIONS
}

\

\centerline{A.R. Mirotin, E. A. Mirotin}

\

\centerline{ amirotin@yandex.ru}  

\

{\small  \textbf{Abstract}.  The purpose of this work is to ascertain
when  arithmetic operations with periodic functions whose domains
may not coincide with the whole real line preserve periodicity.

\medskip

 \textbf{1. Introduction and preliminaries.} The problem under
research is when the arithmetic operations with periodic functions
of one real variable whose domains may not coincide with the  real
line will give periodic functions. The answer is well known in the
case when two nonconstant periodic functions are defined and
continuous on the whole real line and the operation is addition.
In this case the sum is periodic if and only if the periods of
summands are commensurable. But it may be false if the domains of
summands are proper subsets of  reals.

In the following the function $f$ defined on the set $D\subset \mathbb{R}$ is called  \textit{ periodic}  (or  \textit{$T$-periodic}) if $\quad
D+T=D$ and $f(x+T)=f(x)\mbox{ for all } x\in D$ hold for some real
number $T\not=0$. In this case $D$ is called $T$\textit{-invariant}
(or $T$\textit{-periodic}), and $T$ is  called a  \textit{period} of $f$
and $D$. The periods will be always assumed to be positive unless
otherwise stated. The smallest positive period of $f$ and $D$ (if
such exists) is called \textit{fundamental}.

If $D$ is $T$-invariant and   $f(x+T)=f(x)\mbox{ for a. e.}\quad  x\in D$
only, we say that $f$ is  \textit{a. e. periodic}  (with period $T$).

The function  $f$ with domain $D(f)$ is called  \textit{not a.e.
constant} if for every $c$ the set $\{x\in D(f)|f(x)\not=c\}$ has
a positive Lebesgue measure.  $\mathbb{N}, \mathbb{Z}, \mathbb{Q}$, and
$\mathbb{R}$ stand for sets of natural numbers,  integers, rational
numbers, and reals respectively, $\mu$ stands for Lebesgue measure
on $\mathbb{R}$. Below we require that all functions considered have
nonempty domains.

For our purposes the question on the commensurability of  periods of periodic function
is important.  The following example shows that the answer to
this question may be negative.
\par
\medskip
   \textbf{ Example 1.}   Let  $D:=\mathbb{Z}+\sqrt{2}\mathbb{Z}+\sqrt{3}\mathbb{Z}$.  The
function on $D$ defined by the equality
$$
f(k+l \sqrt{2}+m\sqrt{3})=(-1)^m
$$
\noindent
is  bounded and has two incommensurable periods $1$ and $\sqrt{2}$.
\par
    At the same time, the next statement is well known (see, e. g., [3]).

\medskip
   \textbf{ Theorem A.} \textit{If a periodic function $f$ is continuous and nonconstant
on $D(f)$, then $f$ has fundamental period. In particular, every
two periods of$f$ are commensurable.
}
\par
We mention two other conditions, which are  sufficient for
commensurability of the periods of  periodic function.

\par
\medskip
\textbf{Theorem B.} \textit{ Consider a set $D\ne \mathbb{R},\  intD\ne
\emptyset$. If $D$ is periodic then it has the fundamental period.
In particular, if  a periodic function $f$ is defined on $D$, then
$f$ has the fundamental period, too. Thus every two periods of $f$
are commensurable.}

\par
  \textit{Proof}. The set $G$ of all periods
of $D$ is an additive subgroup of
 $\mathbb{R}$ (we consider negative periods and zero as a period of  $D$, too).
Suppose that $G$ is not discrete. Then it is dense in $\mathbb{R}$
(see e.g. [1]). Choose $a \notin D$. The set $-intD+a$  intersects
$G$, so  $-d+a=t$ for some $d$ in $D$ and $t$ in $G$, i.e. $a=d+t
\in D$, a contradiction. Therefore $G= T_0\mathbb{Z}$ for some
$T_0\in \mathbb{R}, \ T_0\ne 0$ (see ibid). This completes the proof.
\par
\medskip
    \textbf{Theorem C.}  \textit{If an a. e. periodic function $f$ is defined, measurable,
and not a.e. constant  on a set $D$ of positive measure, then every two
periods of $f$ are commensurable.}

\par

     \textit{Proof}. Let $T_1,T_2$ be two periods of  $f , S=\mathbb{R}(mod T_1 )$,
i.e. $S$ is a circle with radius $r=T_1/2\pi$. The set $D_1=D(mod
T_1)$ is the subset of $S$ of positive measure. One can assume
that $f$ is defined on $D_1$. The rotation $R_O^\alpha$ of $S$
with the angle $\alpha=T_2/r=(T_2/T_1)2\pi$, which maps $D_1$ on
itself, corresponds to the shift $x\mapsto x+T_2$ of the real
line. If $T_1$ and $T_2$ are incommensurable, $R_O^\alpha$ is an
ergodic transformation of $D_1$ by virtue of the equation
$\alpha/2\pi=T_2/T_1$ (see, e.g., [7], Section II. 5). Since the
function  $f$ on $D_1$ is $R_O^\alpha$-invariant (that is
$f(R_O^\alpha x)=f(x)$ for a. e. $x\in D_1$), it is an a. e.
constant  ([7], ibid), a contradiction.
\par

Note that Burtin's Theorem [2], [4] could be used to prove Theorem
C, too.

\bigskip
\textbf{ 2. Sums of several periodic functions with the common domain.}
    It is well known that the sum of two continuous periodic functions on
$\mathbb{R}$ is periodic if and only if their periods are commensurable.
In this section, we study the periodicity of sums of several periodic
functions $f_i (i=1,...,n)$ in the case where $D(f_1)=\ldots =D(f_n)$
may not coincide
with $\mathbb{R}$. The following example shows that the situation in this
case is more complicated.
\par
\medskip
   \textbf{ Example 2.} Let $D:=\mathbb{Z}+\sqrt{2}\mathbb{Z}+\sqrt{3}\mathbb{Z}$  as in
Example 1. Two functions on $D$ defined by the equalities
$$
f_1(k+l\sqrt{2}+m\sqrt{3})=\frac{1}{|l|+1}-\frac{1}{|m|+1},\quad
f_2(k+l\sqrt{2}+m\sqrt{3})=\frac{1}{|k|+1}+\frac{1}{|m|+1}
$$
\noindent
are bounded and periodic, their periods are incommensurable, but the sum
$f_1+f_2$ is periodic.
\par

    If the periods $T_i$ of several periodic functions $f_i \ (i=1,...,n)$ are
commensurable, it is easy to prove that the sum  $f_1+\cdots +f_n$ is periodic. The converse
is false, in general. If, say, $f_1+f_2= const$, the sum  $f_1+f_2+f_3$  is periodic
for incommensurable $T_1$ and $T_3$. So for converse we should assume that all the
sums of  $f_i$'s where the number of summands is less than $n$ are nonconstant.
\par

\medskip
   \textbf{ Theorem 1. }\textit{ Let $f_1 , f_2 ,\ldots , f_n$ be continuous periodic
functions, which are nonconstant    on their common domain $D$.
If all the sums of  $f_i$'s where the number of summands is less than $n$ are
nonconstant, then the sum $f_1 +\ldots  + f_n$ is periodic if and only if  the
periods of the summands  are commensurable.}
\par
    \textit{Proof}. We shall prove this theorem by induction with the following
additional statement: in the case when the sum is nonconstant the
periods of the summands are commensurable with the period of the
sum. First we shall prove the conclusion of the theorem for $n=2$.
\par
    Suppose that $T_1, T_2$, and $T$ are periods of $f_1, f_2$,
and $f_1+f_2$respectively. Then we have for all $x\in D$
$$
f_1(x+T)+ f_2(x+T)=f_1(x)+ f_2(x),
$$
 \noindent
or
$$
f_1(x+T) -  f_1(x) = f_2(x) -   f_2(x+T) . \eqno (1)
$$

\par
a) Suppose that both sides in (1)
 are nonconstant. Since the
left-hand  side and the right-hand one in (1)
 have
periods $T_1$ and  $T_2$  respectively, Theorem A implies that
 these periods are
commensurable.  Further since $T_1$ and  $T_2$ are commensurable,
the sum  $f_1+f_2$  has  certain period $T^*$ which is
commensurable with $T_1$ and $T_2$.  If $f_1+f_2$ is nonconstant,
then  $T$ and $T^*$ are commensurable by Theorem A, too.
\par
    b) Assume that both sides in (1)
equal to a nonzero constant $c$. The  iteration of the equation
$$
f_1(x+T)-f_1(x)=c \eqno(2)
$$

\noindent implies $f_1(x+nT+mT_1)=f_1(x)+nc$ for all
$m,n\in\mathbb{Z}$.We can find integers $n_k$ and $m_k$, with
$n_k\to\infty$ such that $x+n_kT+m_kT_1\to x$ and we have a
contradiction with the continuity of $f_1$ if $c\not=0$.
\par
c) If both sides of (1) are  zero, then $f_i(x+T)=f_i(x)$, and
$T_i$  and $T$  are commensurable by  Theorem A($i=1,2$).
\par
Now, let the conclusion of the theorem be true for all integers
between 2 and $n$. We shall prove it for $n+1$. Two cases are
possible:
\par
1) The sum $f_1+\cdots +f_{n+1}$ is constant. Then $f_1(x)+\cdots +f_{n+1}(x)=c$
and $f_1(x)+\cdots +f_n(x)=c-f_{n+1}(x)$. Because the left-hand side is
nonconstant, the inductive hypothesis implies that the periods of
$f_1,\ldots ,f_n$  and $T_{n+1}$ are pairwise commensurable.
\par
2) This sum is nonconstant and $T$-periodic. If $g_i(x):=f_i(x+T)-f_i(x)$, then
$$
g_1(x)+\cdots +g_{n+1}(x)=0. \eqno(3)
$$
\noindent
If some $g_i$ is a constant, then it equals to $0$ by b).
\par
  2.1) Let $g_i$'s be nonconstant for $i=1,..., n. $
\par
    d) If the sum $g_1+\cdots +g_n (=-g_{n+1})$ has not proper subsums which
are constant, the periods $T_1, \ldots ,T_n$ are commensurable by inductive
hypothesis. Then the first summand of the sum $(f_1+\cdots +f_n)+f_{n+1}$ has
period of the form  $mT_1$, and again by inductive hypothesis $T_1,T_{n+1}$,
and $T$ are commensurable.
\par
    e) If the sum $g_1+\cdots + g_n (=-g_{n+1})$ has proper subsums which are
constant, let us choose a minimal one, say, $g_1+\cdots +g_k =
const\ (k>1)$. Then by inductive hypothesis, $T_1,\ldots ,T_k$ are
commensurable. Like in d) the first summand of the sum
$(f_1+\cdots +f_k)+(f_{k+1}+\cdots +f_{n+1})$ has the period of
the form  $mT_1$, and by inductive hypothesis $T_1, T_{k+1},
\ldots ,T_{n+1}$ and $T$ are commensurable.

\par
    2.2) If there exist constants among $g_i$'s (which are equal to 0), then
let us renumber  the functions such that $g_1,\ldots ,g_k\ne 0$ and
$g_{k+1}= \ldots =g_{n+1}=0$  where  $k<n+1$. Since for $i$
between $k+1$ and $n+1$ the difference $f_i(x+T)-f_i(x)$ equals to
0, then by Theorem A  the numbers $T_i$ and $T$ are commensurable.
In addition we have  $f_1(x+T)+\cdots +f_k(x+T)=f_1(x)+\cdots
+f_k(x)$ where $k<n+1$. By the hypothesis of the theorem this sum
is nonconstant, so by the inductive hypothesis the periods
$T_1,\ldots ,T_k$ are commensurable with $T$. Moreover, as we have
shown numbers $T_{k+1},\ldots ,T_{n+1}$ are commensurable with
$T$, too.
\par
    We will employ the following lemma to prove Theorem 2. (As was mentioned
    by the referee, one can  prove Theorem 2 using the Proposition 1 in [5] (see also
    [6]);    we give an independent proof which seems to be more elementary).
\par
\medskip
\textbf{Lemma 1.}  \textit{ Let the function $\psi$
 be measurable on the segment $I$. There is a sequence
$\xi_k\downarrow 0$ such that for every sequence $\delta_k,
\delta_k\in (0,\xi_k)$
$$
\lim\limits_{k\to\infty}\psi(x+\delta_k)= \psi(x)   \eqno(4)
$$
for a. e. $x\in I$.}

For the proof see, e. g., [8], proof of Theorem 1.4, especially
formula (1.18).

\medskip
   \textbf{ Theorem 2.}\textit{ Let a. e. $T_i$-periodic functions
$f_i\quad(i=1,\ldots ,n)$ be defined,  measurable, and not a. e.
constant   on the measurable set  $D$ of positive measure. Suppose
that all the sums of $f_i$'s where the number of summands is less
than $n$ are not a. e. constant. The sum $f_1 +\cdots  + f_n$ is
a. e. periodic if and only if the periods of the summands are
commensurable.}

     \textit{Proof}. As in proof of Theorem 1 we shall prove this theorem by induction
with the following additional statement: in the case when the sum
is not a. e. constant the periods of the summands are
commensurable with the period of the sum. First we shall prove the
conclusion of the theorem for $n=2$.
\par
Suppose that $T_1, T_2$, and $T$ are  periods of $f_1,f_2$ and
$f_1+f_2$  respectively.  Then (1) holds for a. e. $x\in D$.

\par
 a) Suppose that both sides in (1) are
 not a. e. constant. Since the left-hand side
and the right-hand one in (1) have  periods $T_1$ and $T_2$
respectively,  Theorem C implies that these periods  are
commensurable. Further the sum  $f_1+f_2$ is defined on the set of
 positive measure. Since $T_1$
and $T_2$  are commensurable, the sum has certain  period $T^*$
which is commensurable with $T_1$ and $T_2$. If $f_1+f_2$ is not
a. e. constant, then $T$ and $T^*$ are commensurable by Theorem C,
too.
\par
    b) Suppose that both sides in (1) equal to a constant $c$ a. e., so that (2)
holds for a. e. $x\in D \quad\ (T$ is the period of $f_1+f_2$).
Then $D$ is $T$-invariant and $T_1$-invariant. Let
$\psi(x)=f_1(x)$ for $x\in D$ and $\psi(x)=0$ for $x\in
\mathbb{R}\setminus D$. We have $\mu(D\cap I)>0$ for some segment
$I\subset \mathbb{R}$. Let $\xi_k\downarrow 0, \xi_k<T_1$ be as in
Lemma 1. If $T$ and $T_1$ are incommensurable one can choose
sequences $m_k, n_k \in \mathbb{Z}$ with the property
$\delta_k:=n_kT+m_kT_1\in(0,\xi_k)$. Then $n_k\ne 0$. Choose $x\in
D$ which satisfies the following three conditions: (4) holds, (2)
holds for $y=x+iT+jT_1$ instead of $x$ for arbitrary integers $i,
j$, and $f_1(y+T_1)=f_1(y)$ for the same $y$. Then $x+\delta_k\in
D$ and the equation (4) implies that
$$
\lim\limits_{k\to\infty}f_1(x+n_kT+m_kT_1)=f_1(x).\eqno(5)
$$
On the other hand, (2) implies that for all $k$
$$
f_1(x+n_kT+m_kT_1)=f_1(x)+n_kc.
$$
\noindent
It follows that $c=0$ and therefore $f_1(x+T) =f_1(x)$ for a. e. $x\in D
$.  Now Theorem C implies that $T$ and $T_1$  are
commensurable. A contradiction. The same is true for $T_2$.

\par
Now, let the conclusion of the theorem be true  for all integers
between $2$ and $n$. We shall  prove it for $n+1$. Two cases are
possible:
\par
1) The sum $f_1+\cdots +f_{n+1}$ is a. e. constant. Then
$f_1(x)+\cdots +f_{n+1}(x)=c$ and $f_1(x)+\cdots
+f_n(x)=c-f_{n+1}(x)$ a. e.. So, by the inductive hypothesis the
periods of  $f_1,\ldots ,f_n$ and the period $T_{n+1}$ of their
sum  are pairwise commensurable (their sum is not a. e. constant
by the hypothesis of the theorem).
\par
    2) This sum is not a. e. constant. Let $g_i(x):=f_i(x+T)-f_i(x)$. Then
$$
g_1(x)+\cdots +g_{n+1}(x)=0.
$$
\par
    2.1) Let $g_i$'s be not a. e. constant for $i=1,..., n+1.$
\par
    c) If the sum $g_1+\cdots +g_n (=-g_{n+1})$ has not proper subsums which
are a. e. constant, the periods $T_1, \ldots ,T_{n+1}$ are
commensurable by inductive hypothesis. Then the sum $f_1+\cdots
+f_{n+1}$ has period of the form  $mT_1$, and again by inductive
hypothesis $T_1$, and $T$ are commensurable.
\par
    d) If the sum $g_1+\cdots + g_n (=-g_{n+1})$ has proper subsums which
are a. e. constant, let us choose a minimal one, say, $g_1+\cdots
+g_k = const$  a. e.  $(k>1)$. Then by inductive hypothesis,
$T_1,\ldots ,T_k$ are commensurable. The first summand of the sum
$(f_1+\cdots +f_k)+(f_{k+1}+\cdots +f_{n+1})$ has the period of
the form  $mT_1$, and by inductive hypothesis $T_1, T_{k+1},
\ldots ,T_{n+1}$ and $T$ are commensurable.
\par
    2.2) If there exist a. e. constants among $g_i$'s for $i=1, ..., n+1,$ say
$g_1= c$ a. e., like in b) it follows that $c=0$ and  $T_1$ is commensurable with
$T$ by Theorem C. So $mT_1 = lT$. Since the sum
$$
f_2+\cdots +f_{n+1}=\sum\limits_{i=1}^{n+1} f_i-f_1
$$

\noindent
is $lT$-periodic and not a. e. constant by inductive hypothesis,
$T_2,\ldots ,T_{n+1}$ , and  $T$
are commensurable by inductive hypothesis, too.
\bigskip

\bigskip
\textbf{ 3. The product  of two periodic functions with possibly
different domains.}  In this section, we assume, as usual, that
the product (and the sum) of several functions with possibly
different domains is defined on the intersection of the domains.
   First consider the following

\textbf{Example 3.} Let $D_1:=\mathbb{Z}+\sqrt{2}\mathbb{Z}+\sqrt{3}\mathbb{Z}+
\sqrt{5}\mathbb{Z}, \ D_2:=\mathbb{Z}+\sqrt{2}\mathbb{Z}+ \sqrt{3}\mathbb{Z}+
\sqrt{7}\mathbb{Z}$. The function $g_1$ on $D_1$ defined by the
equality
$$
g_1(k+ l\sqrt{2}+ m\sqrt{3}+n\sqrt{5})=(|k|+1)(|m|+1)
$$
\noindent has periods $a\sqrt{2} + b\sqrt{5} \ (a,b\in \mathbb{Z})$,
and the function $g_2$ on $D_2$ defined by the equality
$$
g_2(k+ l\sqrt{2}+ m\sqrt{3}+n\sqrt{7})=(|l|+1)/(|m|+1)
$$
\noindent has periods $a  + b\sqrt{7} \ (a,b\in \mathbb{Z})$. But the
product $g_1g_2$ is defined on the set  $D_1\cap
D_2=\mathbb{Z}+\sqrt{2}\mathbb{Z}+ \sqrt{3}\mathbb{Z}$ and has period
$\sqrt{3}$.

At the same time for  $D_i$ with nonempty interior there is a
positive result.

\medskip
    \textbf{Theorem 3.}  \textit{Let $g_i$  be  continuous $T_i$-periodic
functions, and the restrictions  $g_i|intD(g_i)\ne const\
(i=1,2)$. The product $g_1g_2$  is periodic if and only if the
periods $T_1$ and $T_2$ are commensurable.}

We need several lemmas to prove the theorem.

    \medskip
    \textbf{Lemma 2.} \textit{Let $f_i$ be $T_i$-periodic continuous function $(i=1,..,n),
\quad\sum_{i=1}^n f_i\ne  const,\quad D\subseteq\cap_{i=1}^n
D(f_i)$. If the restriction $\sum_{i=1}^n f_i |D$ is $T$-periodic,
then the numbers $T_1^{-1}, \ldots ,T_n^{-1}$, and $T^{-1}$ are
linearly dependent over $\mathbb{Q}.$}
\par
     \textit{Proof.} Assume on the contrary that numbers $T_1^{-1}, \ldots ,
T_n^{-1}$, and $T^{-1}$ are linearly independent over $\mathbb{Q}$.
Since $T/T_1, \ldots , T/T_n$ and  1  are linearly independent
over $\mathbb{Q}$, too, Kronecker Theorem (see e.g. [1], Chapter 7,
section 1, Corollary 2 of Proposition 7) implies, that for $x$ in
$D$ , for every $y$ in $\cap_{i=1}^n D(f_i)$ and $k$ in $\mathbb{N}$
there exist such numbers $q_k$ and $p_{ik}$ in $\mathbb{Z}$, that
$$
|q_kT/T_i - p_{ik} - (y - x)/T_i| < 1/(k \max T_i ) \quad (i = 1, ... ,n)
$$
\noindent
and so
$$
|q_kT - p_{ik} T_i - (y - x)| < 1/k \quad (i = 1, ... ,n).
$$
\noindent
Therefore
$$
\lim_{k\to\infty} (q_k T - p_{ik} T_i ) = y - x  \quad(i = 1,... ,n).
$$
\noindent
Because for $x$ in $D$
$$
f_1(x + q_kT - p_{1k} T_1)  + \cdots  + f_n(x + q_kT - p_{nk} T_n) = f_1 (x) +
\cdots  + f_n(x)
$$
\noindent
and $f_i$'s are continuous, it follows that
$$
 f_1(y) + \cdots + f_n(y) = f_1(x) + \cdots  + f_n(x)
$$
\noindent
and so $\sum_{i=1}^n f_i  = const$. A contradiction.
\par
\medskip

   \textbf{ Corollary 1.}  \textit{Let $f$ be nonconstant continuous $T_1$-periodic
function on $D(f)$. If its restriction to  a subset $D$ of $D(f)$
is $T$-periodic,  then $T$ and $T_1$ are commensurable.}
\par

\medskip
\textbf{Lemma 3.}  \textit{If the set $D_1\ne \mathbb{R}$
 is $T_1$-invariant and its
subset $D, intD\ne\emptyset$, is $T$-invariant, then $T$ and $T_1$
are commensurable.}
\par

     \textit{Proof}. Let us suppose the contrary. Then the set $G = T_1\mathbb{Z} +
T\mathbb{Z}$ is dense in $\mathbb{R}$ by Dirichlet Theorem. Note that
every shift by the element of $G$ maps $D$ into $D_1$. Choose
$a\notin D_1$. Since the open set $a-intD$ intersects $G$,
$a-d=t$,where  $d\in D, t\in G$. Then $a=d+t$ belongs to
 $D_1$, a
contradiction.

\medskip
The following lemma is of intrinsic interest.

\textbf{Lemma 4.}  \textit{Let $f_i$ be $T_i$-periodic  nonconstant continuous
functions with  open domains $D_i\quad(i=1,2)$. The sum $f_1+f_2$
 is periodic
if and only if   the periods of $f_i$' s are commensurable.}
 \par
  \textit{  Proof.} In view of Theorem 1
and Lemma 3 it remains to consider the  case
$D_1\ne\mathbb{R},D_2=\mathbb{R}$. Let $T$ be the period  of the sum
$f_1+f_2$, and suppose that  $T$ and $T_2$ are incommensurable. By
Lemma 3 $mT=kT_1$ for some $m, k$ from $\mathbb{Z}$. Replacing $mT$
by $kT_1$ in the first summand of the left-hand side of the
equality
$$
f_1(x+mT) + f_2(x+mT) = f_1(x) + f_2(x), x\in D_1
$$
\noindent
we have
  $$
f_2(x+mT) = f_2(x), x\in D_1 .
$$
\noindent
It follows from Corollary 1 that  $T_2$ and $T$ are commensurable, a
contradiction.

\textit{Proof of Theorem} 3. First note that the restrictions
$g_i|intD(g_i)$ are $T_i$-periodic, too. So we can assume that
$D(g_i)$ are open. Then the sets
$$
D_i:=\{x\in D(g_i)|g_i(x)\ne0\}\quad (i=1,2)
$$
\noindent are open and $T_i$-invariant. Several  cases are
possible.
\par
1) $D_1\cap D_2\ne\emptyset$. Since $g_1g_2$ is  periodic, the
function on  $D_1\cap D_2$
$$
\log{|g_1g_2|}=\log{|g_1|}+\log{|g_2|}
$$
\noindent
is periodic, too.
\par
1.1). Let both functions $|g_i|$ be nonconstant. Then their periods are
commensurable by Lemma 4.
\par
    1.2). Let both functions $|g_i|$ be constants. Then $D_i\ne\mathbb{R}$
for $i=1,2$ and one can use Lemma  3.
\par
1.3). Let $|g_1|$ is nonconstant, and $|g_2|$  is constant (and so
$g_2(x)=\pm c\ne 0$ ). It was noted above that $D_2\ne\mathbb{R}$. In
view of Lemma 3 we may assume that $D_1= \mathbb{R}$, i.e. $g_1(x)$
has a fixed sign. Let $T$ be the period of $g_1g_2$, so that for
all $x$ in $D(g_2)$ we have
$$
g_1(x+T)g_2(x+T)=g_1(x)g_2(x).\eqno(6)
$$
\noindent Thus the numbers $g_2(x+T)$  and $g_2(x)$ have the same
sign, too, and therefore coincides. Now $T_2$ and $T$ are
commensurable  by  Theorem B.  Then the equality (6) implies
$g_1(x + T) = g_1(x)$ for all  $x$  in $D(g_2)$, and the numbers
 $T_1$
and $T$ are commensurable by  Corollary 1.

\par
2) $D_1\cap D_2=\emptyset$. Suppose that $T_1$  and $T_2$ are
incommensurable. Then for $d_2\in D_2$ one can find two integers
$m,n$  such that $mT_1+nT_2\in  D_1 - d_2$. Therefore
$d_2+nT_2=d_1+(-m)T_1$ for some $d_1\in D_1$. This is impossible
because the left-hand side of the last equality belongs to $D_2$,
but the right-hand one belongs to $D_1$. This completes the proof.
\par
\medskip

    \textbf{Corollary 2.}  \textit{Let $g_i$  be  continuous $T_i$-periodic
functions, and the restrictions  $g_i|intD(g_i)\ne const\
(i=1,2)$. The quotient $g_1 /g_2$  is periodic if and only if the
periods  $T_1$ and  $T_2$ are commensurable.
}

\medskip
\textbf{Remarks.} 1) Let $f_i$ be periodic functions defined on the open
subsets  $D_i\subseteq\mathbb{R}, D_1\ne \mathbb{R}$ and $E_i$ the range
of $f_i (i=1,...,n)$. If the function $F(y_1,\ldots ,y_n)$ on
$E_1\times\ldots\times E_n$ "really depends" of each $y_i$,  the
composition $F(f_1(x),\ldots ,f_n(x))$ is periodic if and only if
the periods of  $f_i$'s are commensurable.  It follows from Lemma
3 immediately. 

For $n=1$ we have the following simple

\textbf{Proposition 1.} \textit{Let the function $F$ on
$\mathbb{R}$ be $T$-periodic.}

1) \textit{If 
$$
f(x+L)-f(x)=nT\quad\forall x\in D(f)\eqno(7)
$$
for some constants $L\ne 0$ and $n\in \mathbb{Z},$ then the composition $F\circ f$ is $L$-periodic.
}

2) \textit{Conversely, if  the composition $F\circ f$ is $L$-periodic, the restriction $F|[0,T)$
 is injective, and $f$ is continuous
on $\mathbb{R},$ then (7) holds for some constants $L\ne 0$ and $n\in \mathbb{Z}.$}

\textit{Proof.} 1) It is obvious.

2) Let $F(f(x+L))=F(f(x))\ \forall x\in \mathbb{R}.$  Since $F|[0,T)$
 is injective, it follows that $f(x+L)-f(x)=n(x)T$ for some function $n:\mathbb{R}\to \mathbb{Z}.$
Now the continuity of $f$ implies that $n=\mathrm{const},$ which completes the proof.

In general the problem on the periodicity of the
composition seems to be open.

2) The problems of  generalization of Theorem 3 for $n>2$
multipliers, for discontinuous multipliers and for general
$D(g_i)$ seem to be open, too.

\par
\medskip
Results of this work were published in \cite{9}.

\end{document}